\documentclass[12pt]{article}

\usepackage{amsthm,amsmath,amssymb}
\usepackage{url}
\usepackage{srcltx} 
\usepackage{enumerate}
\usepackage[utf8]{inputenc}
\usepackage{mathrsfs}
\usepackage{xspace}
\usepackage{framed}
\usepackage{enumerate}
\usepackage{xcolor}
\usepackage{xspace}
\usepackage{footnote}
\makesavenoteenv{tabular}
\usepackage[colorlinks=true,citecolor=black,linkcolor=black,urlcolor=blue]{hyperref}
\usepackage{cleveref}

\def\abs#1{\left| #1 \right|}
\def\paren#1{\left( #1 \right)}
\def\acc#1{\left\{ #1 \right\}}

\renewcommand{\le}{\leqslant}
\renewcommand{\ge}{\geqslant}

\theoremstyle{plain}
\newtheorem{theorem}{Theorem}

\theoremstyle{definition}

\newtheorem{conjecture}[theorem]{Conjecture}

\theoremstyle{remark}
\newtheorem{remark}[theorem]{Remark}



\begin{document}

\title{Doubled patterns with reversal and square-free doubled patterns}

\author{
Antoine Domenech\\
\small ENS de Lyon, Lyon, France\\[-0.8ex]
\small\tt antoine.domenech@ens-lyon.fr\\
\and
Pascal Ochem\\
\small LIRMM, CNRS, Universit\'e de Montpellier\\[-0.8ex]
\small Montpellier, France \\[-0.8ex] 
\small\tt ochem@lirmm.fr
}
\maketitle

\begin{abstract}
In combinatorics on words, a word $w$ over an alphabet $\Sigma$ is said to avoid a pattern $p$ over an alphabet $\Delta$
if there is no factor $f$ of $w$ such that $f=h(p)$ where $h:\Delta^*\to\Sigma^*$ is a non-erasing morphism.
A pattern $p$ is said to be $k$-avoidable if there exists an infinite word over a $k$-letter alphabet that avoids $p$.
A pattern is \emph{doubled} if every variable occurs at least twice. Doubled patterns are known to be $3$-avoidable.
Currie, Mol, and Rampersad have considered a generalized notion which allows variable occurrences to be reversed.
That is, $h(V^R)$ is the mirror image of $h(V)$ for every $V\in\Delta$.
We show that doubled patterns with reversal are $3$-avoidable.
We also conjecture that (classical) doubled patterns that do not contain a square are $2$-avoidable.
We confirm this conjecture for patterns with at most 4 variables.
This implies that for every doubled pattern $p$, the growth rate of ternary words avoiding $p$ is at least
the growth rate of ternary square-free words.
A previous version of this paper containing only the first result has been presented at WORDS 2021.
\end{abstract}

\section{Introduction}\label{sec:intro}

The \emph{mirror image} of the word $w=w_1w_2\dots w_n$ is the word $w^R=w_nw_{n-1}\dots w_1$.
A pattern with reversal $p$ is a non-empty word over an alphabet\\ $\Delta=\acc{A,A^R,B,B^R,C,C^R\dots}$
such that $\acc{A,B,C,\dots}$ are the \emph{variables} of $p$.
An \emph{occurrence} of $p$ in a word $w$ is a non-erasing morphism $h:\Delta^*\to\Sigma^*$
satisfying $h(X^R)=(h(X))^R$ for every variable $X$ and such that $h(p)$ is a factor of $w$.
The avoidability index $\lambda(p)$ of a pattern with reversal $p$ is the size of the
smallest alphabet $\Sigma$ such that there exists an infinite word $w$ over $\Sigma$ containing no occurrence of $p$.
A pattern $p$ such that $\lambda(p)\le k$ is said to be $k$-avoidable.
To emphasive that a pattern is without reversal (i.e., it contains no $X^R$), it is said to be \emph{classical}.
A pattern is \emph{doubled} if every variable occurs at least twice.

Our aim is to strengthen the following result.
\begin{theorem}~\cite{BellGoh:2007,Ochem:2004,O16}\label{doubled}
Every doubled pattern is $3$-avoidable.
\end{theorem}

First, we extend it to patterns with reversal.
\begin{theorem}\label{doubledr}
Every doubled pattern with reversal is $3$-avoidable.
\end{theorem}

Then, we notice that all the known classical doubled patterns that are $2$-unavoidable contain a square,
such as $AABB$, $ABAB$, or $ABCCBADD$.

\begin{conjecture}\label{conj_sqf}
Every square-free doubled pattern is $2$-avoidable.
\end{conjecture}

Notice that \Cref{conj_sqf} is related to but independent of the following conjecture.

\begin{conjecture}~\cite{O16,OchemPinlou}\label{conj_finite}
There exist only finitely many 2-unavoidable doubled patterns.
\end{conjecture}

The proof of \Cref{conj_sqf} for patterns up to 3 variables follows from the $2$-avoidability of
$ABACBC$, $ABCBABC$, $ABCACB$ and $ABCBAC$.
We were able to verify it for patterns up to 4 variables.

\begin{theorem}\label{4var}
Every square-free doubled pattern with at most 4 variables is $2$-avoidable.
\end{theorem}

Finally, we obtain a lower bound on the number of ternary words avoiding a doubled pattern. 
The factor complexity of a factorial language $L$ over $\Sigma$ is $f(n)=\abs{L\cap\Sigma^n}$.
The growth rate of $L$ over $\Sigma$ is $\lim_{n\to\infty}f(n)^{\tfrac1n}$.
We denote by $GR_3(p)$ the growth rate of ternary words avoiding the doubled pattern $p$.

\begin{theorem}\label{doublede}
For every doubled pattern $p$, $GR_3(p)\ge GR_3(AA)$.
\end{theorem}

Let $v(p)$ be the number of distinct variables of the pattern~$p$.
In the proof of Theorem~\ref{doubled}, the set of doubled patterns is partitioned as follows:
\begin{enumerate}
 \item Patterns with $v(p)\le3$: the avoidability index of every ternary pattern has been determined~\cite{Ochem:2004}.
 \item Patterns shown to be $3$-avoidable with the so-called power series method:
 \begin{itemize}
  \item Patterns with $v(p)\ge6$~\cite{BellGoh:2007}
  \item Patterns with $v(p)=5$ and prefix $ABC$ or length at least $11$~\cite{O16}
  \item Patterns with $v(p)=4$ and prefix $ABCD$ or length at least $9$~\cite{O16}
 \end{itemize}
 \item Ten sporadic patterns with $4\le v(p)\le5$ whose $3$-avoidability cannot be deduced from the previous results:
 they have been shown to be $2$-avoidable~\cite{O16} using the method in~\cite{Ochem:2004}.
\end{enumerate}

The proof of \Cref{doubledr,doublede} use the same partition.
\Cref{v3,v6,v45} are each is devoted to one type of doubled pattern with reversal.
\Cref{4var} is proved in \Cref{sec_4}
\Cref{doublede} is proved in \Cref{sec_e}

\section{Preliminaries}
A word $w$ is \emph{$d$-directed} if for every factor $f$ of $w$ of length $d$, the word $f^R$ is not a factor of $w$.
\begin{remark}\label{remark}
If a $d$-directed word contains an occurrence $h$ of $X.X^R$ for some variable $X$, then $|h(X)|\le d-1$.
\end{remark}

A variable that appears only once in a pattern is said to be \emph{isolated}.
The \emph{formula} $f$ associated to a pattern $p$ is obtained by replacing every isolated variable in $p$ by a dot.
The factors between the dots are called \emph{fragments}.
An occurrence of a formula $f$ in a word $w$ is a non-erasing morphism $h$ such that the $h$-image of every fragment of $f$ is a factor of $w$.
As for patterns, the avoidability index $\lambda(f)$ of a formula $f$ is the size of the smallest alphabet allowing the existence of an infinite word containing no occurrence of~$f$.
Recently, the avoidability of formulas with reversal has been considered by Currie, Mol, and Rampersad~\cite{CMR2017,CMR2018} and Ochem~\cite{O_rev_infty}.

Recall that a formula is \emph{nice} if every variable occurs at least twice in the same fragment.
In particular, a doubled pattern is a nice formula with exactly one fragment.

The \emph{avoidability exponent} $AE(f)$ of a formula $f$ is the largest real $x$ such that every $x$-free word avoids $f$.
Every nice formula $f$ with $v(f)\ge3$ variables is such that $AE(f)\ge1+\tfrac1{2v(f)-3}$~\cite{OchemRosenfeld2021}.

Let $\simeq$ be the equivalence relation on words defined by $w\simeq w'$ if $w'\in\acc{w,w^R}$.
Avoiding a pattern up to $\simeq$ has been investigated for every binary formulas~\cite{CM2020}.
Remark that for a given classical pattern or formula $p$,
avoiding $p$ up to $\simeq$ implies avoiding simultaneously all the variants of $p$ with reversal.

Recall that a word is $(\beta^+,n)$-free if it contains no repetition with exponent strictly greater than $\beta$ and period at least $n$.

\section{Formulas with at most 3 variables}\label{v3}
For classical doubled patterns with at most 3 variables, all the avoidability indices are known.
There are many such patterns, so it would be tedious to consider all their variants with reversal.

However, we are only interested in their 3-avoidability, which follows from the 3-avoidability
of nice formulas with at most 3 variables~\cite{OchemRosenfeld2017}.

Thus, to obtain the 3-avoidability of doubled patterns with reversal with at most 3 variables,
we show that every minimally nice formula with at most 3 variables is 3-avoidable up to $\simeq$.

The minimally nice formulas with at most 3 variables, up to symmetries, are determined in~\cite{OchemRosenfeld2017} and listed in the following table.
Every such formula $f$ is avoided by the image by a $q$-uniform morphism of either any infinite $\paren{\frac54^+}$-free word $w_5$ over $\Sigma_5$
or any infinite $\paren{\frac75^+}$-free word $w_4$ over $\Sigma_4$, depending on whether the avoidability exponent of $f$ is smaller than $\frac75$.

\begin{center}
\begin{tabular}{|l|l|l|l|l|l|l|}
\hline
Formula $f$ & $=f^R$ & $AE(f)$ & Word & $q$ & $d$ & freeness\\ \hline
$ABA.BAB$ & yes & 1.5 & $g_a(w_4)$ & 9 & 9 & $\paren{\tfrac{131}{90}^+,28}$\\ \hline
$ABCA.BCAB.CABC$ & yes & 1.333333333 & $g_b(w_5)$ & 6 & 8 & $\paren{\tfrac43^+,25}$\\ \hline
$ABCBA.CBABC$ & yes & 1.333333333 & $g_c(w_5)$ & 4 & 9 & $\paren{\tfrac{30}{23}^+,18}$\\ \hline
$ABCA.BCAB.CBC$ & no  & 1.381966011 & $g_d(w_5)$ & 9 & 4 & $\paren{\tfrac{62}{45}^+,37}$\\ \hline \hline
$ABA.BCB.CAC$ & yes & 1.5 & $g_e(w_4)$\footnote{The formula $ABA.BCB.CAC$ seems also avoided up to~$\simeq$ by the Hall-Thue word, i.e.,
the fixed point of $\texttt{0}\to\texttt{012}$; $\texttt{1}\to\texttt{02}$; $\texttt{2}\to\texttt{1}$.} & 9 & 4 & $\paren{\tfrac{67}{45}^+,37}$\\ \hline
$ABCA.BCAB.CBAC$ & yes\footnote{We mistakenly said in~\cite{OchemRosenfeld2017} that $ABCA.BCAB.CBAC$ is different from its reverse.} & 1.333333333 & $g_f(w_5)$ & 6 & 6 & $\paren{\tfrac{31}{24}^+,31}$\\ \hline
$ABCA.BAB.CAC$ & yes & 1.414213562 & $g_g(w_4)$ & 6 & 8 & $\paren{\tfrac{89}{63}^+,61}$\\ \hline
$ABCA.BAB.CBC$ & no & 1.430159709 & $g_h(w_4)$ & 6 & 7 & $\paren{\tfrac{17}{12}^+,61}$\\ \hline
$ABCA.BAB.CBAC$ & no & 1.381966011 & $g_i(w_5)$ & 8 & 7 & $\paren{\tfrac{127}{96}^+,41}$\\ \hline
$ABCBA.CABC$ & no & 1.361103081 & $g_j(w_5)$ & 6 & 8 & $\paren{\tfrac43^+,25}$\\ \hline
$ABCBA.CAC$ & yes & 1.396608253 & $g_k(w_5)$ & 6 & 13 & $\paren{\tfrac43^+,25}$\\ \hline
\end{tabular}\label{formula}
\end{center}
In the table above, the columns indicate respectively, the considered minimally nice formula $f$, whether $f$ is equivalent to its reversed formula, the avoidability exponent of $f$,
the infinite ternary word avoiding $f$, the value $q$ such that the corresponding morphism is $q$-uniform, the value such that the avoiding word is $d$-directed,
the suitable property of $(\beta^+,n)$-freeness used in the proof that $f$ is avoided.
We list below the corresponding morphisms.

\noindent
\begin{minipage}[b]{0.2\linewidth}
\centering
$$\begin{array}{c}
 g_a\\
 \texttt{002112201}\\
 \texttt{001221122}\\
 \texttt{001220112}\\ 
 \texttt{001122012}\\ 
\end{array}$$
\end{minipage}
\begin{minipage}[b]{0.2\linewidth}
\centering
$$\begin{array}{c}
 g_b\\
 \texttt{021221}\\
 \texttt{021121}\\
 \texttt{020001}\\ 
 \texttt{011102}\\ 
 \texttt{010222}\\ 
\end{array}$$
\end{minipage}
\begin{minipage}[b]{0.2\linewidth}
\centering
$$\begin{array}{c}
 g_c\\
 \texttt{2011}\\
 \texttt{1200}\\
 \texttt{1120}\\ 
 \texttt{0222}\\
 \texttt{0012}\\ 
\end{array}$$
\end{minipage}
\begin{minipage}[b]{0.2\linewidth}
\centering
$$\begin{array}{c}
 g_d\\
 \texttt{020112122}\\
 \texttt{020101112}\\
 \texttt{020001222}\\ 
 \texttt{010121222}\\
 \texttt{000111222}\\ 
\end{array}$$
\end{minipage}
\begin{minipage}[b]{0.2\linewidth}
\centering
$$\begin{array}{c}
 g_e\\
 \texttt{001220122}\\
 \texttt{001220112}\\
 \texttt{001120122}\\ 
 \texttt{001120112}\\
\end{array}$$
\end{minipage}
\begin{minipage}[b]{0.16\linewidth}
\centering
$$\begin{array}{c}
 g_f\\
 \texttt{012220}\\
 \texttt{012111}\\
 \texttt{012012}\\ 
 \texttt{011222}\\
 \texttt{010002}\\ 
\end{array}$$
\end{minipage}
\begin{minipage}[b]{0.16\linewidth}
\centering
$$\begin{array}{c}
 g_g\\
 \texttt{021210}\\
 \texttt{011220}\\
 \texttt{002111}\\ 
 \texttt{001222}\\ 
\end{array}$$
\end{minipage}
\begin{minipage}[b]{0.16\linewidth}
\centering
$$\begin{array}{c}
 g_h\\
 \texttt{011120}\\
 \texttt{002211}\\
 \texttt{002121}\\ 
 \texttt{001222}\\ 
\end{array}$$
\end{minipage}
\begin{minipage}[b]{0.16\linewidth}
\centering
$$\begin{array}{c}
 g_i\\
 \texttt{01222112}\\
 \texttt{01112022}\\
 \texttt{01100022}\\ 
 \texttt{01012220}\\
 \texttt{01012120}\\ 
\end{array}$$
\end{minipage}
\begin{minipage}[b]{0.16\linewidth}
\centering
$$\begin{array}{c}
 g_j\\
 \texttt{021121}\\
 \texttt{012222}\\
 \texttt{011220}\\ 
 \texttt{011112}\\
 \texttt{000102}\\ 
\end{array}$$
\end{minipage}
\begin{minipage}[b]{0.16\linewidth}
\centering
$$\begin{array}{c}
 g_k\\
 \texttt{022110}\\
 \texttt{021111}\\
 \texttt{012222}\\ 
 \texttt{012021}\\
 \texttt{011220}\\ 
\end{array}$$
\end{minipage}

\noindent
As an example, we show that $ABCBA.CAC$ is avoided by $g_k(w_5)$.
First, we check that $g_k(w_5)$ is $\paren{\tfrac43^+,25}$-free using the main lemma in~\cite{Ochem:2004}, that is,
we check the $\paren{\tfrac43^+,25}$-freeness of the $g_k$-image of every $\paren{\tfrac54^+}$-free word of length at most $\frac{2\times\tfrac43}{\tfrac43-\tfrac54}=32$.
Then we check that $g_k(w_5)$ is $13$-directed by inspecting the factors of $g_k(w_5)$ of length $13$.
For contradiction, suppose that $g_k(w_5)$ contains an occurrence $h$ of $ABCBA.CAC$ up to~$\simeq$.
Let us write $a=|h(A)|$, $b=|h(B)|$, $c=|h(C)|$.

Suppose that $a\ge25$. Since $g_k(w_5)$ is $13$-directed, all occurrences of $h(A)$ are identical.
Then $h(ABCBA)$ is a repetition with period $|h(ABCB)|\ge25$.
So the $\paren{\tfrac43^+,25}$-freeness implies the bound $\tfrac{2a+2b+c}{a+2b+c}\le\tfrac43$, that is, $a\le b+\tfrac12c$.

In every case, we have $$a\le\max\acc{b+\tfrac12c,24}.$$
Similarly, the factors $h(BCB)$ and $h(CAC)$ imply $$b\le\max\acc{\tfrac12c,24}$$
and $$c\le\max\acc{\tfrac12a,24}.$$
Solving these inequalities gives $a\le36$, $b\le24$, and $c\le24$.
Now we can check exhaustively that $g_k(w_5)$ contains no occurrence up to~$\simeq$ satisfying these bounds.

Except for $ABCBA.CBABC$, the avoidability index of the nice formulas in the above table is 3.
So the results in this section extend their $3$-avoidability up to~$\simeq$.

\section{The power series method}\label{v6}
The so-called power series method has been used~\cite{BellGoh:2007,O16} to prove the $3$-avoidability of many classical doubled patterns
with at least $4$ variables and every doubled pattern with at least $6$ variables, as mentioned in the introduction.

Let $p$ be such a classical doubled pattern and let $p'$ be a doubled pattern with reversal obtained by adding some $-^R$ to $p$.
Witout loss of generality, the leftmost appearance of every variable $X$ of $p$ remains free of $-^R$ in $p'$.
Then we will see that $p'$ is also $3$-avoidable. The power series method is a counting argument that relies on the following observation.
If the $h$-image of the leftmost appearance of the variable $X$ of $p$ is fixed, say $h(X)=w_X$, then there is exactly one possibility
for the $h$-image of the other appearances of $X$, namely $h(X)=w_X$. This observation can be extended to $p'$, since there is also exactly one possibility
for $h(X^R)$, namely $h(X^R)=w_X^R$.

Notice that this straightforward generalization of the power series method from classical doubled patterns to doubled patterns with reversal
cannot be extended to avoiding a doubled pattern up to~$\simeq$. Indeed, if $h(X)=w_X$ for the leftmost appearance of the variable $X$ and $w_X$
is not a palindrome, then there exist two possibilities for the other appearances of $X$, namely $w_X$ and $w_X^R$.

\section{Sporadic patterns}\label{v45}
Up to symmetries, there are ten doubled patterns whose $3$-avoidability cannot be deduced by the previous results.
They have been identified in~\cite{O16} and are listed in the following table.
\begin{center}
\begin{table}\label{sporadic}\caption{The seven sporadic patterns on 4 variables and the three sporadic patterns on 5 variables}
\begin{tabular}{|l|l|}
\hline
Doubled pattern & Avoidability exponent\\ \hline
$ABACBDCD$ & 1.381966011\\ \hline
$ABACDBDC$ & 1.333333333\\ \hline
$ABACDCBD$ & 1.340090632\\ \hline
$ABCADBDC$ & 1.292893219\\ \hline
$ABCADCBD$ & 1.295597743\\ \hline
$ABCADCDB$ & 1.327621756\\ \hline
$ABCBDADC$ & 1.302775638\\ \hline
$ABACBDCEDE$ & 1.366025404\\ \hline
$ABACDBCEDE$ & 1.302775638\\ \hline
$ABACDBDECE$ & 1.320416579\\ \hline
\end{tabular}
\end{table}
\end{center}

Let $w_5$ be any infinite $\paren{\frac54^+}$-free word over $\Sigma_5$ and let $h$ be the following $9$-uniform morphism.
$$
\begin{array}{c}
h(0)=020022221\\
h(1)=011111221\\
h(2)=010202110\\
h(3)=010022112\\
h(4)=000022121\\
\end{array}
$$
First, we check that $h(w_5)$ is $7$-directed and $\paren{\tfrac{139}{108}^+,46}$-free.
Then, using the same method as in \Cref{v3}, we show that $h(w_5)$ avoids up to~$\simeq$ these ten sporadic patterns simultaneously.

\section{Square-free doubled patterns with at most 4 variables}\label{sec_4}
Here we show \Cref{4var}, that is, every square-free doubled pattern with at most 4 variables is 2-avoidable.
We list them as follows:
\begin{itemize}
 \item Among patterns that are equal up to letter permutation, we only list the lexicographically least.
 \item If a pattern is distinct from its mirror image, we only list the lexicographically least among the pattern and its mirror image. 
 \item We do not list patterns that contain a square-free doubled pattern as a strict factor.
 \item We do not list patterns that contain an occurrence of $ABACBC$, $ABCACB$, $ABCBABC$, $ABCDBDABC$, $ABCDBDAC$, $ABACDCBD$, or their mirror image.
 \item We do not include the seven sporadic patterns on 4 variables from Table~1, which are 2-avoidable.
\end{itemize}

Table 2 contains every pattern $p$ in this list with an infinite binary word avoiding $p$.
Let us detail how to read Table 2:
\begin{itemize}
 \item A morphism is $m$ given in the format $m(0)/m(1)/...$ 
 \item We denote by $b_2$, $b_4$, $b_5$ the famous morphisms $01/10$, $01/21/03/23$, $01/23/4/21/0$, respectively.
 \item We denote by $w_k$ any infinite $RT(k)^+$-free word over $\Sigma_k$.
 \item If the avoiding word is a pure morphic word $m^{\omega}(0)$, then $m$ is given.
 \item If the avoiding word is a morphic word $f(m^{\omega}(0))$, then we write $m$; $f$.
 \item If the avoiding word is of the form $f(w_k)$, then we write $w_k$; $f$.
\end{itemize}
The proofs that a (pure) morphic word avoids a pattern use Cassaigne's algorithm~\cite{cassaignealgo}
and the proofs that a morphic image word a Dejean word avoids a pattern use the technique described in \Cref{v3}.

\begin{center}
\begin{table}\tiny\label{123}\caption{Binary words avoiding doubled patterns}
\begin{tabular}{|l|l|l|}
\hline
Doubled pattern & Avoiding word\\ \hline
ABCABDCBD & $w_5$; 0010101110/0010011000/0001111110/0001110101/0000011001\\ \hline
ABCACDCBD & $w_5$; 000101010111/000100110111/000011001111/000001011111/000000111111\\ \hline
ABCBABDBCBD & $b_4$; 01/00/10/11\\ \hline
ABCBADCBCD & $b_4$; 0000/0011/1111/1010\\ \hline
ABCBDACBCD & $b_4$; 01/00/10/11\\ \hline
ABCBDBCACBCD & $b_2$\\ \hline
ABCBDCBACBCD & $b_4$; 1000/0111/0110/0010\\ \hline
ABCDACBD & $w_5$; 00100110111111000/00100110111011000/\\
 & 00011110110101010/00001111111011010/00001010101011111\\ \hline
ABCDBACBD & $w_6$; 010101111100/010010100000/001001110111/000111111101/000101010111/000100011011\\ \hline
ABCDBADC & $w_5$; 10001000101111101010110/00000110110101000111111/\\ 
 & 00000101011100100111111/00000011101010010011111/00000011011000101011111\\ \hline
ABCDBCBACBD & 001/011 \\ \hline
ABCDCACBD & $b_5$; 0011110110000/0011010100110/0001111100111/0001110001000/0001101101111\\ \hline
ABCDCBABCD & avoided by every $\paren{\tfrac32^+,4}$-free binary word~\cite{IOS:2004}\\ \hline
ABCDCBCACBD & $b_5$; 00/01/10/110/111\\ \hline
\hline
ABACDCBCD & $w_5$; 10011011000/01011111000/00111010100/00100100111/00001111111\\ \hline
ABCABDBCD & $w_5$; 0010111111/0010011110/0010011100/0000010101/0000001101\\ \hline
ABCADBCBD & $w_5$; 001011010000/001001111000/000110011001/000011101010/000010111111\\ \hline
ABCADCBCD & $w_5$; 001101111000/001101101000/001001111111/000101110101/000001100101\\ \hline
ABCBADBDC & $w_5$; 0011111110110/0001010111100/0000101101110/0000011010111/0000001011111\\ \hline
ABCBDABCD & $w_4$; 1111/1101/0010/0000\\ \hline
ABCBDABDC & $w_5$; 101110000001/101100100001/011111110100/010001111110/010001101110\\ \hline
ABCBDACBD & $b_5$; 00/01/10/110/111\\ \hline
ABCBDADBC & $w_5$; 00110111010010/00110000000010/00011111111011/00011110101000/00010101100011\\ \hline
ABCBDADBDC & $b_5$; 111/101/000/011/001\\ \hline
ABCBDBABDBC & $b_5$; 00/01/10/110/111\\ \hline
ABCBDBABDC & $b_5$; 000/011/001/111/101\\ \hline
ABCBDBACBCD & $b_4$; 01/00/10/11\\ \hline
ABCBDBACD & $w_5$; 0001111101010/0001110111000/0001011111111/0000111001111/0000011011001\\ \hline
ABCBDBADBDC & 011/100\\ \hline
ABCBDBADC & 00111101110000/00111011000010/00111010100000/00011001001111/00010101111111\\ \hline
ABCBDBCABCD & $b_5$; 00/01/10/110/111\\ \hline
ABCBDBCACBD & $b_5$; 00/01/10/110/111\\ \hline
ABCBDBCAD & $w_5$; 00011110110011/00011101101001/00011011010100/00010111111110/00000011111010\\ \hline
ABCBDBCBABCD & $b_4$; 000/111/10/01\\ \hline
ABCBDBCBACBCD & $b_2$\\ \hline
ABCBDBCBACBD & $b_4$; 00/01/10/11\\ \hline
ABCBDBCBACD & 001/110\\ \hline
ABCBDCABCD & $b_5$; 00/10/111/01/011\\ \hline
ABCBDCABD & $w_5$; 10000000011/01111010010/01101100010/01011111110/00001010101\\ \hline
ABCBDCACBD & $b_5$; 111/101/000/100/110\\ \hline
ABCBDCBABCD & $b_5$; 00/01/10/110/111\\ \hline
ABCBDCBACBD & $b_5$; 00/01/10/110/111\\ \hline
ABCBDCBACD & $b_5$; 00/01/10/1100/111\\ \hline
ABCBDCBAD & $w_5$; 001101101100/001011111111/001001111100/000110010100/000001110100\\ \hline
ABCBDCBCABD & $b_4$; 000/111/10/01\\ \hline
ABCBDCBCAD & $w_5$; 1111100/1100110/0110101/0010010/0000101\\ \hline
ABCDADCB & $w_5$; 0000010001111110101000100111110111/0000010001111100100001100101101111/\\
 & 0000001001111111010000110101111011/0000001001111110110100010101111011/\\
 & 0000000101110010000111111010010111\\ \hline
ABCDBABDC & $w_5$; 0011111110101/0010110111010/0010101110000/0000111111001/0000110110001\\ \hline
ABCDBADBC & $w_5$; 01011111111/01001000111/00101000011/00011110101/00000001011\\ \hline
ABCDBCACBD & $b_5$; 101/000/110/111/100\\ \hline
ABCDBCBACD & $w_5$; 0110101/0100000/0011110/0001111/0000111\\ \hline
ABCDBCBAD & $w_5$; 00010111001010/00001111010101/00001110001010/00001100111111/00001100010110\\ \hline
ABCDBDAC & $w_5$; 00000011011011001110001111011010110000101111010100100101110111/\\
 & 00000011011011000010011110110101000010101111010100100101110111/\\
 & 00000010110011110101010011000111000010101111010100100101110111/\\
 & 00000010101101101000100011111101000010101111010100100101110111/\\
 & 00000010101011001110001111010011000010101111010100100101110111\\ \hline
ABCDBDADBC & $w_5$; 01111101/00111100/00111001/00110110/00000101\\ \hline
ABCDCACDB & $w_5$; 00110001000110/00101011111110/00011111010011/00010101011111/00000001010011\\ \hline
\end{tabular}
\end{table}
\end{center}

\section{Growth rate of ternary words avoiding a doubled pattern}\label{sec_e}
\Cref{doublede} obviously holds for $p=AA$.
Without loss of generality, we do not need to consider a doubled pattern $p$ that contains an occurrence of another doubled pattern.
In particular, $p$ is square-free.
So we need to show that $GR_3(p)$ is at least $GR_3(AA)$, which is close to $1.30176$~\cite{Shur}.

If $p$ is 2-avoidable, then $p$ is avoided by sufficiently many ternary words.
By Lemma 4.1 in~\cite{Ochem:2004}, $\lambda(p)=2$ implies that $GR_3(p)\ge2^{\tfrac12}>GR_3(AA)$.
Thus, \Cref{conj_sqf} implies \Cref{doublede}.
By \Cref{4var}, we can assume that $v(p)\ge5$.
We can also rule out the three sporadic patterns on 5 variables from Table 1, which are 2-avoidable.

According to the partition of the set of doubled patterns mentioned in the introduction,
there remains to consider the doubled patterns $p$ whose 3-avoidability has been obtained via the power series method.
In that case, we even get $GR_3(p)>2>GR_3(AA)$.

%
%
\section{Conclusion}\label{conclusion}

Unlike classical formulas, we know that there exist avoidable formulas with reversal of arbitrarily high avoidability index~\cite{O_rev_infty}.
Maybe doubled patterns and nice formulas are easier to avoid. We propose the following open problems.

\begin{itemize}
 \item Are there infinitely many doubled patterns up to~$\simeq$ that are not $2$-avoidable?
 \item Is there a nice formula up to~$\simeq$ that is not $3$-avoidable?
\end{itemize}

A first step would be to improve \Cref{doubledr} by generalizing the $3$-avoidability of doubled patterns with reversal to doubled patterns up to~$\simeq$.
Notice that the results in \Cref{v3,v45} already consider avoidability up to~$\simeq$.
However, the power series method gives weaker results.
Classical doubled patterns with at least $6$ variables are $3$-avoidable because 
$$1-3x+\paren{\frac{3x^2}{1-3x^2}}^v$$
has a positive real root for $v\ge6$.
The (basic) power series for doubled patterns up to~$\simeq$ with $v$ variables would be 
$$1-3x+\paren{\frac{6x^2}{1-3x^2}-\frac{3x^2+3x^4}{1-3x^4}}^v.$$
The term $\frac{6x^2}{1-3x^2}$ counts for twice the term $\frac{3x^2}{1-3x^2}$ in the classical setting, for $h(V)$ and $h(V)^R$.
The term $\frac{3x^2+3x^4}{1-3x^4}$ corrects for the case of palindromic $h(V)$, which should not be counted twice.
This power series has a positive real root only for $v\ge10$.
This leaves many doubled patterns up to~$\simeq$ whose $3$-avoidability must be proved with morphisms.

Looking at the proof of \Cref{doubledr}, we may wonder if a doubled pattern with reversal is always easier to avoid than the corresponding classical pattern.
This is not the case: backtracking shows that $\lambda(ABCA^RC^RB)=3$, whereas $\lambda(ABCACB)=2$~\cite{Ochem:2004}.

To get a more precise version of both conjectures~\ref{conj_sqf} and~\ref{conj_finite}, we plan to obtain the (conjectured) list of all 2-unavoidable doubled patterns,
which should be a finite list containing no square-free pattern.

\end{document}